\input ./style/arxiv-general.cfg
\documentclass[MSNbibl,citesort,number,seceqn,dvips]{arxbj}
\makeatletter
   \@ifpackageloaded{graphicx}{}{\usepackage{graphicx}}
\makeatother
\usepackage{multirow,dcolumn}
\usepackage{graphicx}
\usepackage{stfloats}

\aid{0}
\volume{21}
\issue{2}
\pubyear{2015}
\firstpage{1238}
\lastpage{1259}
\doi{10.3150/14-BEJ604} 

\makeatletter
  \fnbelowfloat
\newcommand{\rrvert}{\vert}
\newcommand{\llvert}{\vert}
\newtheorem{them}{Theorem}[section]
\newproclaim{assu}{Assumption}[section]
\newtheorem{prop}{Proposition}[section]
\newtheorem{Le}{Lemma}[section]

\newcommand{\bdelta}{{\boldsymbol{\delta}}}

\newcommand{\bphi}{{\boldsymbol{\phi}}}

\newcommand{\btheta}{{\boldsymbol{\theta}}}

\newcommand{\bxi}{{\boldsymbol{\xi}}}

\newcommand{\bu}{{\mathbf{u}}}

\newcommand{\bv}{{\mathbf{v}}}

\newcommand{\xh}{\hat{X}}

\newcommand{\oh}{1/2}

\newcommand{\R}{\mathbb{R}} 

\newcommand{\g}{\gamma}

\newcommand{\eps}{\varepsilon}

\newcommand{\epsh}{\hat{\varepsilon}}
\newcommand{\E}{\mathrm{E}}

\newcommand{\mtoinf}{m\rightarrow\infty}

\newcommand{\eqref}[1]{(\ref{#1})}

\newcommand{\ARMA}{\operatorname{ARMA}}
\newcommand{\AR}{\operatorname{AR}}
\newcommand{\Laplace}{\operatorname{Laplace}}
\newcommand{\MA}{\operatorname{MA}}
\newcolumntype{d}[1]{D{.}{.}{#1}}
\makeatother

\begin{document}
\begin{frontmatter}

\title{Reaction times of monitoring schemes for ARMA time series}
\runtitle{Monitoring ARMA time series}

\begin{aug}
\author[A]{\inits{A.}\fnms{Alexander} \snm{Aue}\corref{}\thanksref{A,e1}\ead[label=e1,mark]{aaue@ucdavis.edu}},
\author[A]{\inits{C.}\fnms{Christopher} \snm{Dienes}\thanksref{A,e2}\ead[label=e2,mark]{crdienes@ucdavis.edu}},\\
\author[B]{\inits{S.}\fnms{Stefan} \snm{Fremdt}\thanksref{B}\ead[label=e3]{stefan.fremdt@ruhr-uni-bochum.de}}
\and
\author[C]{\inits{J.}\fnms{Josef} \snm{Steinebach}\thanksref{C}\ead[label=e4]{jost@math.uni-koeln.de}}
\address[A]{Department of Statistics, University of California, One
Shields Avenue, Davis, CA 95616, USA.\\
\printead{e1,e2}}

\address[B]{Department of Mathematics, Institute of Statistics, Ruhr-University Bochum, D-44780 Bochum, Germany.
\printead{e3}}

\address[C]{Mathematisches Institut, Universit\"at zu K\"oln, Weyertal
86-90, D-50931 K\"oln, Germany.\\
\printead{e4}}
\end{aug}

\received{\smonth{4} \syear{2013}}
\revised{\smonth{7} \syear{2013}}

%
\begin{abstract}
This paper is concerned with deriving the limit distributions of
stopping times devised to sequentially uncover structural breaks in the
parameters of an autoregressive moving average, ARMA, time series. The
stopping rules are defined as the first time lag for which detectors,
based on CUSUMs and Page's CUSUMs for residuals, exceed the value of a
prescribed threshold function. It is shown that the limit distributions
crucially depend on a drift term induced by the underlying ARMA
parameters. The precise form of the asymptotic is determined by an
interplay between the location of the break point and the size of the
change implied by the drift. The theoretical results are accompanied by
a simulation study and applications to electroencephalography, EEG, and
IBM data. The empirical results indicate a satisfactory behavior in
finite samples.
\end{abstract}

%

\end{frontmatter}

\section{Introduction}\label{sec1}

Sequential change-point analysis is concerned with uncovering in an
on-line fashion what is called structural breaks, deviations from a
pre-specified in-control scenario. For the case of time series,
relevant for this paper, the natural in-control scenario is the
stationarity of the underlying stochastic process. More traditionally,
sequential change-point techniques were developed for breaks in the
mean and variance in sequences of independent observations. The
corresponding literature is reviewed in the monographs Basseville and Nikiforov \cite{bn93} and Cs\"org\H{o} and Horv\'ath \cite{ch97}. A more recent
survey of both sequential and historical procedures is given in
Aue and Horv{\'a}th \cite{ah12+}.

The particular approach to sequential change-point analysis of this
paper is grounded in the work of Chu \textit{et al.} \cite{csw96}, who developed
procedures using a training sample to estimate an initial model and to
monitor for deviations from that model as soon as new observations
arrive. This contribution, originally written with applications to
econometric data in mind, has been extended in a number of ways.
Further sequential procedures covering financial time series were
discussed in Andreou and Ghysels \cite{ag06}, and Aue \textit{et al.} \cite{ahr09}. Berkes
\textit{et al.} \cite{bghk04} introduced methodology applicable to GARCH processes.
Gombay and Serban \cite{gs09} worked with autoregressive processes, while
Gombay and Horv{\'a}th \cite{gh09} considered weakly stationary time series.
Refinements using bootstrap were considered in Kirch \cite
{k08} and Hu\v
{s}kov\'a and Kirch \cite{hk11}, while resampling schemes were studied by
Hu\v{s}kov\'a \textit{et al.} \cite{hkps08}.

The basic time series model being utilized in this paper is the class
of linear autoregressive moving average, ARMA, processes made popular
through the works of Box \textit{et al.} \cite{bjr08}. ARMA processes find widespread
applications in a number of fields as evidenced, for example, in the
recent text Shumway and Stoffer \cite{ss10}. As advocated by Brown \textit{et al.}
\cite{bde75}) in a regression setting, the proposed monitoring procedures are
based on the residuals obtained from an ARMA model fit to the original
data based on a training sample of size $m$ for which stationarity of
the underlying process is assumed. If the process remains stationary
after the monitoring starts, then residuals of the training period and
the monitoring period should possess similar properties. The test
procedures to be introduced here are based on traditional cumulative
sum, CUSUM, statistics and a modification, Page's CUSUM statistics (see
Page \cite{p54,p55}). The latter tend to react faster to deviations from
the in-control scenario and satisfy certain optimality criteria (see Lorden
\cite{l71}).
CUSUMs for residuals of ARMA processes were discussed in a
retrospective setting in Bai \cite{b93}, Yu \cite{y07} and Robbins \textit{et al.}
\cite{rgla11}, and in a sequential framework in Dienes and Aue \cite{da12+}. Recent
work on Page's CUSUMs can be found in Fremdt \cite{f12+a,f12+b}.

A stopping rule is then defined as a first crossing time, that is, the
time lag for which either the CUSUM or Page's CUSUM statistic exceed a
threshold value tolerable for the in-control case. The focus of this
paper is on deriving the asymptotic distributions of these stopping
rules for the situation that deviations from stationarity of the
underlying process occur. The particular deviations of interest are the
classic change in mean and general changes in the second-order
dynamics, with an emphasis on changes in the variance (or scale) due to
the nature of the data examples provided in this paper. Namely, the
finite-sample properties of the proposed methods are discussed in two
case studies. The first of the applications involves EEG data. Here
interest is in detecting the occurrence of an epileptic seizure (see
Davis \textit{et al.} \cite{dlr06}). The second application deals with closing prices
of IBM stock, a classic data set that has been analyzed with historical
procedures for the presence of breaks in variance (see, e.g.,
Tsay \cite{t88}). Accompanying simulation evidence indicates that the
procedure works satisfactory for these two examples.

The paper is organized as follows. Section~\ref{sec2} details the ARMA
model and states the hypotheses to be tested. Section~\ref{sec3}
quantifies the large-sample behavior of the delay times incurred by the
CUSUM and Page's CUSUM procedure. Applications to EEG and IBM data are
discussed in Section~\ref{sec:apps}. All proofs are given in
Section~\ref{sec5}.

\section{The model}\label{sec2}

Let $\mathbb{Z}$ denote the set of integers. In what follows,
$(Y_t\dvtx  t\in\mathbb{Z})$ denotes the $\ARMA(p,q$) process
specified by the stochastic recurrence equations
%
\begin{eqnarray}\label{ARMA}
\phi_t(B) (Y_t - \mu_t) =
\theta_t(B)\eps_t, \qquad t\in\mathbb{Z},
\end{eqnarray}
where $\mu_t$ are mean parameters, $\phi_t(z)= 1-\phi_{t,1} z -
\cdots- \phi_{t,p} z^p$ and $\theta_t(z) = 1+\theta_{t,1} z +\cdots
+ \theta_{t,q} z^q$ denote respectively the autoregressive and moving
average polynomials, and $B$ the backshift operator. The innovations
$(\eps_t\dvtx  t\in\mathbb{Z})$ are assumed to be independent random
variables with zero mean and variance $\sigma_t^2$. As usual, it is
further required that $\phi_t$ and $\theta_t$ have no common zeroes
and that the ARMA process is causal and invertible, which means
%
\begin{eqnarray}\label{inv-cau}
\phi_t(z)\neq0 \quad \mbox{and}\quad  \theta_t(z)\neq0\qquad \mbox{for
all }|z|\leq1.
\end{eqnarray}
The parametric model in \eqref{ARMA} depends on the parameter vectors
$\bxi_t = (\mu_t,\bphi_t,\btheta_t,\sigma_t)^\prime$, where
$\bphi_t = (\phi_{t,1},\ldots,\phi_{t,p})^\prime$ and $\btheta_t
= (\theta_{t,1},\ldots,\theta_{t,q})^\prime$, with $^\prime$
denoting transposition. These vectors may be time dependent and
interest is in monitoring the constancy of the $\bxi_t$ in a
sequential fashion. This is important because constancy of the $\bxi
_t$ would imply stationarity of the underlying ARMA process, so that
standard methods are available for estimation and prediction purposes.
To set up the monitoring, a training period of size $m+p$ is utilized
for which
%
\begin{eqnarray}\label{NCA}
Y_{1-p},\ldots,Y_m \mbox{ are governed by }
\bxi_t=\bxi_0=(\mu _0,\bphi_0,
\btheta_0,\sigma_0)^\prime.
\end{eqnarray}
As Chu \textit{et al.} \cite{csw96} elaborate, this training period may be used to
estimate the parameters of an initial non-contaminated model and to
express limit results in the form $m\to\infty$. In particular, let
$(X_t\dvtx  t\in\mathbb{Z})$ be the centered sequence defined by $X_t
= Y_t - \mu_t$ and define $\hat{\bxi}_m = (\hat{\mu}_{m},\hat
{\boldsymbol{\phi}}_{m},\hat{\boldsymbol{\theta}}_{m}
,\hat{\sigma}_{m})^\prime$ to be a $\sqrt{m}$-consistent estimator
for $\bxi
_0$ obtained from the training period data. This gives the model residuals
\begin{eqnarray*}
\epsh_t = \xh_t - \sum_{j = 1}^p
\hat{\phi}_{m,j}\xh_{t-j} - \sum_{j =
1}^q
\hat{\theta}_{m,j}\epsh_{t-j},
\end{eqnarray*}
with $\xh_t=Y_t-\hat{\mu}_{m}$ and initializations $\epsh_{-q+1} =
\cdots=
\epsh_0 = 0$ in case $q>0$.

In the following, two sets of hypotheses will be considered. First, the
focus will be on the arguably most studied case for which only mean
breaks are permitted. The sequential testing problem then becomes
\begin{eqnarray*}
&H_0\dvtx & Y_{m+1},Y_{m+2},\ldots\mbox{ have mean }
\mu_0;
\\
&H_A^\mu\dvtx & Y_{m+1},\ldots,Y_{m+k^*-1}\mbox{
have mean }\mu_0, \mbox { but }
\\
&& Y_{m+k^*},Y_{m+k^*+1},\ldots\mbox{ have mean }\mu_A
\neq\mu_0,
\end{eqnarray*}
where here the constancy of the remaining model ARMA parameters is
required, so that changes may only affect the mean. It may sometimes be
of greater importance to test for changes in the underlying
second-order dynamics. This can be done via testing the general
sequential hypotheses
\begin{eqnarray*}
&H_0\dvtx & Y_{m+1},Y_{m+2},\ldots\mbox{ are governed
by } {\bxi}_0;
\\
&H_A^{\bxi}\dvtx & Y_{m+1},\ldots,Y_{m+k^*-1}\mbox{
are governed by } {\bxi }_0,
\mbox{ but }
\\
&& Y_{m+k^*},Y_{m+k^*+1},\ldots\mbox{ are governed by } {
\bxi}_A\neq {\bxi}_0.
\end{eqnarray*}
Under $H_A^{\bxi}$ the decomposition ${\bxi}_A = {\bxi}_0 + \bdelta
_m^{\bxi}$ will be utilized, where $\bdelta_m^{\bxi} = (\delta
_m^\mu,\bdelta_m^\phi,\bdelta_m^\theta,\delta_m^\sigma)^\prime$
denotes the difference in parameter values.

For both sets of hypotheses, one can now proceed as follows. If the
respective null scenarios hold, then the residuals $\epsh_t$ should
roughly resemble the corresponding innovations $\eps_t$ and suitably
constructed statistics should therefore behave similarly on the
training period and after monitoring commences. Under the alternatives,
this should not be the case. This approach will be detailed in the next section.

\section{Monitoring schemes and their large-sample properties}\label{sec3}
\setcounter{equation}{0}

\subsection{CUSUM and Page's CUSUM procedures under the null}\label{sec:3.1}

Testing procedures for the set of hypotheses introduced in the previous
section are commonly defined as stopping times that reject the null if
a detector crosses the boundary prescribed by a threshold function.
Popular choices for the detector are based on cumulative sum, CUSUM,
statistics and on its variant, called Page's CUSUM. Let $\mathbb{N}$
denote the positive integers. To introduce the CUSUM of (squared)
residual procedures, define for $k\in\mathbb{N}$ the detectors
%
\begin{eqnarray}\label{Dmu}
\hat{D}_\mu(m,k) = \sum_{t = m+1}^{m+k}
\epsh_t - \frac{k}{m}\sum_{t =
1}^{m}
\epsh_t  \quad \mbox{and}\quad  \hat{D}_{\bxi}(m,k) = \sum
_{t = m+1}^{m+k}\epsh_t^2
- \frac{k}m \sum_{t = 1}^m
\epsh_t^2.
\end{eqnarray}
The detector $\hat{D}_\mu(m,k)$ is built from the residuals $\epsh
_t$ and used to test $H_0$ against $H_A^\mu$, while the detector $\hat
{D}_\bxi$ is built from the squared residuals $\epsh_t^2$ and used to
test $H_0$ against $H_A^\bxi$. Using the class of weight functions
%
\begin{equation}\label{g-def}
g_\gamma(m,k)= \sqrt{m} \biggl(1+\frac{k}m \biggr) \biggl(
\frac
{k}{m+k} \biggr)^\gamma,
\end{equation}
indexed by a sensitivity parameter $\gamma\in[0,1/2)$, a stopping
time corresponding to the detector $\hat{D}_\mu(m,k)$ can be defined by
%
\begin{eqnarray}\label{tau_m}
\tau_\mu(m) = \min \bigl\{k\in\mathbb{N}\dvtx \bigl|\hat{D}_\mu(m,k)\bigr|
\geq c_\alpha\hat\sigma_m g_\gamma(m,k) \bigr\},
\end{eqnarray}
where $c_\alpha=c_{\alpha}(\gamma)$ is a critical constant, derived
from the limit distribution of the detector under $H_0$ (see Theorem~\ref{Th1} below), ensuring that $P(\tau_\mu(m)<\infty) = \alpha$
for a given level $\alpha\in(0,1)$.

The stopping time $\tau_\bxi(m)$ for the detector $\hat{D}_\bxi
(m,k)$ is defined analogously: Let $\hat{\eta}_m^2$ denote a weakly
consistent estimator of the quantity $\eta^2 = \E[(\eps_1^2 - \sigma
^2)^2]$. Then $\tau_\bxi(m)$ is given by replacing $\hat D_\mu(m,k)$
and $\hat\sigma_m$ with $\hat D_\bxi(m,k)$ and $\hat\eta_m$, respectively.

Page's CUSUM procedure is a modification of the CUSUM detectors in
\eqref{Dmu} based on the adjusted detectors
%
\begin{eqnarray}\label{Dmu^P}
\hat{D}_\mu^{P}(m,k) &=& \max_{0\leq k^\prime\leq k} \bigl|\hat
{D}_\mu(m,k) - \hat {D}_\mu\bigl(m,k^\prime\bigr)
\bigr|,
\nonumber
\\[-8pt]\\[-8pt]
\hat{D}_{\bxi}^{P}(m,k) &=& \max
_{0\leq k^\prime\leq k} \bigl|\hat {D}_{\bxi}(m,k) - \hat {D}_{\bxi}
\bigl(m,k^\prime\bigr) \bigr|,\nonumber
\end{eqnarray}
setting $\hat D_\mu(m,0)=\hat D_\bxi(m,0)=0$. Utilizing the same
class of weight functions in \eqref{g-def} as before gives rise to the
Page-type stopping time
%
\begin{eqnarray}\label{tau_m^p}
\tau_\mu^{P}(m) = \min\bigl\{k\in\mathbb{N}\dvtx
\hat{D}_\mu ^{P}(m,k) \geq c_\alpha^{P}
\hat\sigma_mg_\gamma(m,k)\bigr\},
\end{eqnarray}
where $c_\alpha^{P}=c_\alpha^{P}(\gamma)$ controls again the level
of the sequential procedure. The stopping time $\tau^P_{\bxi}(m)$ is
defined in a similar fashion. These sequential detectors were
introduced in the seminal papers (Page \cite{p54,p55}).

All procedures are based on residuals instead of directly on the
observations. This has the advantage that the notoriously difficult
estimation of long-run variances of the dependent observations can be
completely avoided. Better size and power properties are expected from
this approach as pointed out in Robbins \textit{et al.} \cite{rgla11}, who confirmed
these statements in an extensive simulation study.

The large-sample behavior under the null hypotheses for the four
detectors is quantified in the following two theorems, the first one of
which states the results for the mean only procedures.

\begin{them}\label{Th1}
Let $(Y_t\dvtx  t\in\mathbb{Z})$ follow the ARMA equations \eqref
{ARMA} and assume that $\E[|\eps_1|^\nu]<\infty$ for some $\nu>2$.
Then it holds under $H_0$ and for all real $c$ that
\begin{enumerate}[(b)]
\item[(a)]
$\displaystyle \lim_{\mtoinf} P \biggl({\displaystyle \frac{1}{\hat{\sigma}_m}\displaystyle \sup
_{k\geq
1}\displaystyle \frac{|\hat{D}_\mu(m,k)|}{g_\gamma(m,k)}\leq c} \biggr) = P \biggl(\displaystyle \sup
_{0< x< 1}\displaystyle \frac{\llvert W(x)\rrvert }{x^\g}\leq c \biggr)$,\vspace*{2pt}

\item[(b)]
$\displaystyle \lim_{\mtoinf} P \biggl({\displaystyle \frac{1}{\hat{\sigma}_m}
\displaystyle \sup_{k\geq
1}\displaystyle \frac{\hat{D}_\mu^{P}(m,k)}{g_\gamma(m,k)}\leq c} \biggr) = P \biggl(\displaystyle \sup
_{0< x< 1}\displaystyle \sup_{0\leq y\leq x} \displaystyle \frac{1}{x^\g}\biggl
\llvert W(x)-\displaystyle \frac{1-x}{1-y}W(y)\biggr\rrvert \leq c \biggr)$,
\end{enumerate}
where $(W(x)\dvtx  x\in[0,1])$ denotes a standard Brownian motion.
\end{them}

%
\begin{them}\label{Th2}
Let $(Y_t\dvtx  t\in\mathbb{Z})$ follow the ARMA equations \eqref
{ARMA} and assume that $\E[|\eps_1|^\nu]<\infty$ for some $\nu>4$.
Then, under $H_0$ and for all real $c$, the limit results of Theorem~\ref{Th1} are retained if $\hat D_\mu(m,k)$, $\hat D_\mu^P(m,k)$ and
$\hat\sigma_m$ are replaced with the respective objects $\hat D_\bxi
(m,k)$, $\hat D_\bxi^P(m,k)$ and $\hat\eta_m$.
\end{them}

The proofs of the theorems follow from the results in Dienes and Aue
\cite{da12+} for the CUSUM procedure, and from a combination of the latter
with the proofs in Fremdt \cite{f12+a} for Page's CUSUM
procedure. Tables
containing simulated critical values for a selection of sensitivity
parameters $\g$ and test levels $\alpha$ can be found in Horv\'ath \textit{et
al.} \cite{hhks04} for the limit in Theorem~\ref{Th1}, part (a) and in Fremdt
\cite{f12+a} for the limit in part (b).

\subsection{Limiting delay times for mean breaks}\label{ssec1}

The quality of monitoring procedures is often quantified via the mean
delay time which measures how long, on average, one has to wait before
the structural break in the underlying processes is detected. For
example, certain optimality criteria for Page's CUSUM were developed in
Lorden \cite{l71}. The monograph by Basseville and Nikiforov
\cite{bn93} gives
an account of the subsequent contributions in this area. The main
theoretical contribution of this paper is the derivation of the
complete limit distribution of the stopping times under consideration.
Taking the mean of this distribution, one obtains in particular also
the information on the average delay time. Related results in the
literature are Aue and Horv{\'a}th \cite{ah04}, Aue \textit{et al.} \cite{ahr09} and
Fremdt \cite{f12+a}. To account for the ARMA time series
character, modifications of
the methodology in these papers become necessary. These will be
developed in the following.

It is subsequently assumed that $H_A^\mu$ holds and that thus changes
in the second-order structure of the ARMA process do not occur. Notice
that assumption \eqref{inv-cau} implies that the reciprocals of $\phi
_t(z)$ and $\theta_t(z)$ admit, for $|z|\leq1$, the power series expansions
%
\begin{equation}\label{recip}
\frac{1}{\phi_t(z)} = \sum_{\ell=0}^\infty
\pi_\ell(\bphi _t)z^\ell \quad \mbox{and}\quad
\frac{1}{\theta_t(z)} = \sum_{\ell=0}^\infty
\psi_\ell(\btheta _t)z^\ell.
\end{equation}
Denoting the training period estimates of the autoregressive and moving
average polynomials by $\hat{\phi}_{m}(z)$ and $\hat{\theta
}_{m}(z)$, for large
enough $m$, one finds analogously power series expansions for their
reciprocals. These will be written as
%
\begin{equation}\label{recip-hat}
\frac{1}{\hat{\phi}_{m}(z)} = \sum_{\ell=0}^\infty
\pi_\ell(\hat {\boldsymbol{\phi}}_{m})z^\ell \quad \mbox{and}\quad  \frac{1}{\hat{\theta}_{m}(z)} = \sum_{\ell=0}^\infty
\psi_\ell (\hat{\boldsymbol{\theta}}_{m} )z^\ell.
\end{equation}
Under $H_A^\mu$, the asymptotic behavior of the delay time will depend
on the size of the mean change $\delta^\mu_m=\mu_A - \mu_0$ which
in turn induces the drift term
%
\begin{equation}
\Delta_m^\mu = \delta^\mu_m
\Biggl(1-\sum_{j=1}^p\phi_{0,j}
\Biggr)\sum_{\ell
=0}^\infty\psi_\ell(
\btheta_0) = \delta^\mu_m\frac{\phi_0(1)}{\theta_0(1)}.
\end{equation}
Note that the difference of pre-mean and post-mean is allowed to depend
on $m$, so that one could more explicitly write $\mu_{A,m}$. The
precise limit distribution will crucially depend on the interplay
between the behavior of the drift term $\Delta_m^\mu$ and the
location of
the mean change ${k^{*}}$. This leads to the following set of assumptions
which, in view of the theorems to come, are formulated for a general
sequence $\Delta_m$ and not directly for $\Delta_m^\mu$.
Superscripts, such as
$\mu$ here, will indicate which drift term is being used.

\begin{assu}\label{As1}
It is required that
\begin{itemize}[(b)]
\item[(a)] there is $\theta>0$ such that ${k^{*}}= \lfloor\theta
m^\beta\rfloor$ with $\beta\in[0,1)$, where $\lfloor\cdot\rfloor
$ denotes integer part; 
\item[(b)] $\sqrt{m}|\Delta_m|\to\infty$; 
\item[(c)] $|\Delta_m| = \mathrm{O} ({1} )$. 
\end{itemize}
\end{assu}

Part (a) of Assumption~\ref{As1} specifies the order of the
change-point ${k^{*}}$ as a power of $m$. It is a standard assumption in
the change-point literature. However, it should be noted that the
expression $k^* = \lfloor\theta m^\beta\rfloor$ is not unique for
fixed $m$ and $k^*$, and different specifications of $\theta$ and
$\beta$ may lead to different limit distributions. A discussion of
this matter can be found in Section~3 of Fremdt \cite
{f12+a}. Note also that
parts (b) and (c) implicitly allow for the decay of the sequence
$|\Delta_m
|$. The proofs show that the form of the limit distribution of the
stopping times depends then on the asymptotic behavior of the sequence
$|\Delta_m| m^{\g-1/2}{k^{*}}^{1-\g}$ of scaled drift terms. Due to
part (a)
of Assumption~\ref{As1} which allows for the re-expression of $k^*$ in
terms of $m$, they depend consequently on the asymptotic behavior of
the scaled terms
\[
\tilde\Delta_m=|\Delta_m| m^{\beta(1-\g)-\oh+\g},
\]
which do not explicitly contain ${k^{*}}$ anymore. We distinguish between
the three cases
\[
\mbox{(i) } \tilde\Delta_m\to0,\qquad  \mbox{(ii) } \tilde\Delta_m\to\tilde
C_1\in(0,\infty),\qquad  \mbox{(iii) } \tilde\Delta_m\to\infty.
\]
%
%
In case (ii), it follows from part (a) of Assumption~\ref
{As1} that
$|\Delta_m| {m^{\g-\oh}}{k^{*}}^{1-\g}\to\theta^{1-\g}\tilde
{C}_1 =
C_1\in(0,\infty)$. 
For this scenario and any real $c$ define $d_1 = d_1(c)$ to be the
unique solution of
%
\begin{equation}\label{d1}
d_1 = 1 - \frac{c}{C_1}d_1^{1-\g}.
\end{equation}

In order to exhibit the asymptotic distribution of the stopping times,
introduce first the case-dependent distribution function $\Psi$ by
setting, for all real arguments $u$,
\[
\Psi(u) = \cases{ \Phi(u), &\quad \mbox{in case (i)},
\cr
P \Bigl({\displaystyle \sup
_{d_1 < x<1}W(x)\leq u} \Bigr),&\quad \mbox{in case (ii)},
\cr
P
\Bigl({\displaystyle \sup_{0<x<1}W(x)\leq u} \Bigr) = \cases{0,&\quad $u<0$,
\cr
2
\Phi(u) - 1,& \quad $u\geq0$,}
&\quad \mbox{in case (iii)},}
\]
where $\Phi$ denotes the standard normal distribution function. The
next theorem gives the large-sample behavior of $\tau_\mu(m)$ and
$\tau_\mu^P(m)$.

\begin{them}\label{AD}
Let $(Y_t\dvtx  t\in\mathbb{Z})$ follow the ARMA equations \eqref
{ARMA} so that \eqref{inv-cau} and \eqref{NCA} hold, and suppose that
Assumption~\ref{As1} is satisfied for $\Delta_m=\Delta_m^\mu$.
Then it holds under $H_A^\mu$ for all real $u$ that
\begin{eqnarray*}
\mbox{\textup{(a)} } \lim_{m\rightarrow\infty}P \biggl(\frac{\tau_\mu^P(m) -
a_m(c_\alpha^P)}{b_m
(c_\alpha^P)}\leq u \biggr) = 1-{\Psi}(-u). 
\end{eqnarray*}
Additionally, 
\begin{eqnarray*}
\mbox{\textup{(b)} }\lim_{m\rightarrow\infty}P \biggl(\frac{\tau_\mu(m) -
a_m(c_\alpha)}{b_m
(c_\alpha)}\leq u \biggr) = \Phi(u),
\end{eqnarray*}
where $a_m(c)$ is the unique positive solution of
%
\begin{eqnarray}\label{amc}
a_m(c) = \biggl({\frac{cm^{\oh-\g}}{|\Delta_m^\mu|}+\frac
{{k^{*}}}{(a_m (c))^\g}}
\biggr)^{1/(1-\g)}
\end{eqnarray}
and
\begin{eqnarray*}
b_m(c) = \frac{\sigma\sqrt{a_m(c)}}{|\Delta_m^\mu|} \biggl({1-\g \biggl({1-
\frac{{k^{*}}}{a_m(c)}} \biggr)} \biggr)^{-1}.
\end{eqnarray*}
\end{them}

The proof of Theorem~\ref{AD} is in Section~\ref{Pssec1}. Note that
the uniqueness of $a_m(c)$ follows from a rewriting of equation \eqref
{amc} to
\[
a_m(c)=\frac{cm^{1/2-\gamma}}{|\Delta_m^\mu|}\bigl(a_m(c)
\bigr)^\gamma+k^*.
\]
Now it can be seen that $a_m(c)$ solves an equation of the form $x=
ax^\gamma+ b$ for appropriately chosen $a>0$, $b>0$ and $\gamma\in
[0,1/2)$. Since $a_m(c)>0$, it is unique as the intersection of the
identity with a transformed power function whose exponent is smaller
than one.

A similar result can be obtained for the squared-residual procedures
$\tau_\bxi(m)$ and $\tau_\bxi^P(m)$ after appropriate modification.
The proof of the following theorem may also be found in Section~\ref
{Pssec1} below.

\begin{them}
\label{ADmu^2}
Let $(Y_t\dvtx  t\in\mathbb{Z})$ follow the ARMA equations \eqref
{ARMA} so that \eqref{inv-cau} and \eqref{NCA} hold, and suppose that
Assumption~\ref{As1} is satisfied for $\Delta_m=(\Delta_m^\mu)^2$. Then,
under $H_A^\mu$ for all real $u$, the limit results of Theorem~\ref
{AD} are retained if $\tau_\mu(m)$, $\tau_\mu^P(m)$ and $\sigma$
are replaced with the respective objects $\tau_\bxi(m)$, $\tau_\bxi
^P(m)$ and $\eta$.
\end{them}

Some discussion is in order. First, the limit distributions for Page's
CUSUM and the traditional CUSUM coincide for the early change scenario
(i). Therefore, all procedures work similar in a large-sample
setting. The critical values for the traditional CUSUM are somewhat
smaller than those for Page's CUSUM (comparing the tables in Horv\'ath
\textit{et al.} \cite{hhks04} with those of Fremdt \cite{f12+a}), giving it a
slight edge
for this case. However, limiting distributions are different for the
intermediate and late change scenarios (ii) and (iii), respectively. Here, Page's CUSUM outperforms the traditional
CUSUM. This can be explained by the fact that, unlike Page's CUSUM, the
traditional CUSUM is not resetting and so becomes less sensitive to a
change the later it occurs after the onset of monitoring.\looseness=-1

Second, in view of the last paragraph, Page's CUSUM is generally
preferred for applications unless the changes happen early. For the
early change scenario (i) both procedures perform alike in
finite samples (based on simulations not reported in the paper), but as
the theoretical results indicate, the performance of the traditional
CUSUM decays noticeably for (ii) and (iii). In
fact, this stopping rule often exhibits significant non-zero
probabilities of non-detection in intermediate and late changes
scenarios if the monitoring period is not sufficiently long.

Third, the sensitivity of the test can be adjusted by the statistician
through the choice of $\gamma$. For example, it has been pointed out
by Aue and Horv{\'a}th \cite{ah04} that the term $a_m(c)$ can be interpreted
as the average delay time $E[\tau]$, where $\tau$ stands for any of
the stopping times under consideration. For the early change scenario
(i), it follows then that $E[\tau]\approx(c/|\Delta
_m|)^{1/(1-\gamma)}m^{1-2\gamma/[2(1-\gamma)]}$. This quantity
becomes small if $\gamma$ is chosen close to $1/2$, thus ensuring a
quicker detection. However, there is an obvious trade-off between
detection time and false alarm rates, with the latter increasing with
increasing $\gamma$. Similar computations can be obtained for cases
(ii) and (iii) as well.

\subsection{Limiting delay times for scale breaks}\label{ssec2}

In view of the applications, for which only changes in the scale are
considered, presentation in this section is focused on the case of a
break in the scale parameter $\sigma$ only. All other parameters are
assumed to remain the same before and after the change occurs. The
section closes with remarks for the general case, but a more in-depth
analysis is beyond the scope of the present paper. The special case of
the general alternative $H_A^\bxi$, for which only the scale parameter
is subject to change, will be called $H_A^\sigma$ in the following. A
change of scale will induce the drift term
\[
\Delta_m^\sigma=\bigl(\delta_m^\sigma
\bigr)^2+2\sigma_0\delta_m^\sigma
\]
into the squared-residual procedures. If this drift term satisfies the
regularity conditions imposed through Assumption~\ref{As1}, then the
asymptotic delay time distribution can be quantified accordingly.

\begin{them}
\label{AD-sigma}
Let $(Y_t\dvtx  t\in\mathbb{Z})$ follow the ARMA equations \eqref
{ARMA} so that \eqref{inv-cau} and \eqref{NCA} hold, and suppose that
Assumption~\ref{As1} is satisfied for $\Delta_m= \Delta_m^\sigma
\neq0$ and
$\delta_m^\sigma=\mathrm{O}(1)$. Then the results of Theorem~\ref
{ADmu^2} remain valid under $H_A^{\sigma}$.
\end{them}

The proof of Theorem~\ref{AD-sigma} is given in Section~\ref{Pssec2}.
The general case is much more difficult to handle. The induced drift
term will be a complicated function of the pre-break parameters $\bxi
_0$ and post-break parameters $\bxi_A$. In principle, the arguments
developed in order to verify the theorems of Sections~\ref{ssec1} and
\ref{ssec2} could be adjusted to this case. However, one has to keep
track of additional terms, the number of which may be growing
exponentially in the number of parameters. Given the complexity of the
proofs, we refrain from pursuing this direction further for this paper.

\section{Applications}
\label{sec:apps}

In order to demonstrate the proposed methodology in the finite sample
setting, two case studies are provided in this section. The first
involves an EEG data set considered in Davis \textit{et al.} \cite{dlr06}, the second
is a classic data set on IBM stock given in Box \textit{et al.} \cite{bjr08},
previously analyzed for breaks in variance with retrospective methods.

\subsection{EEG data}

\begin{figure}[b]

\includegraphics{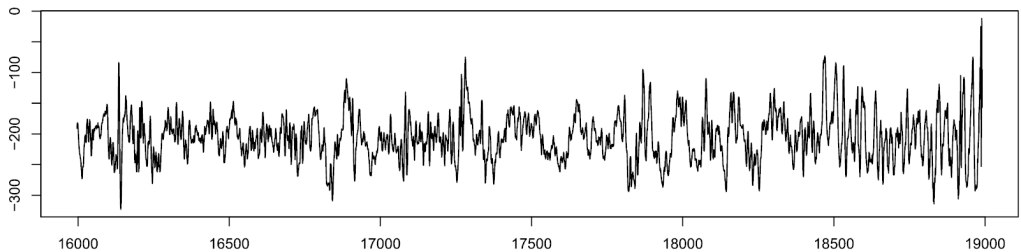}

\caption{EEG data set.}\label{fig:1}
\end{figure}

In this section, the proposed methodology is applied to two snapshots
of a longer series of 32\,768 EEG measurements observed from a female
patient diagnosed with left temporal lobe epilepsy.\footnote{We thank
Dr. Beth Malow (formerly Department of Neurology, University of
Michigan) for providing the data.} This is the ``T3 channel'' data of
Davis \textit{et al.} \cite{dlr06}. Measurements were taken at a sampling rate of 100
Hz (i.e., 100 observations per second), so that the recording took
place over a time period of 5 minutes and 28 seconds. As explained in
Davis \textit{et al.} \cite{dlr06}, expert analysis suggests the onset of an epileptic
seizure at observation 18\,500. Their (retrospective) segmentation
procedure estimates the seizure onset at observation 18\,580. A
similar analysis is reported in Ombao \textit{et al.} \cite{orvsm01}.
Particular interest here is in two segments of the original data
focusing on the interval 16\,000--19\,000 before and immediately
after the suspected seizure onset. These observations are plotted in
Figure~\ref{fig:1}. A visual inspection of the time series plot
indicates that the level of the observations remains roughly the same.
There is, however, an apparent increase in the amplitude around time
18\,500, perhaps indicating a scale break.
%

To test for this possibility, the following two scenarios for the
training period, both of size $m=1000$, were considered:
\begin{enumerate}[(TP2)]
\item[(TP1)]  Observations 16\,001--17\,000,
\item[(TP2)]  Observations 17\,001--18\,000.
\end{enumerate}
The training periods predate the epileptic seizure, with (TP1) implying
a longer monitoring period before the break occurrence than (TP2). The
choices of training periods enable to examine the effect of the
change-point location on the monitoring procedures. In each case, model
selection procedures suggest nearly identical $\AR(4)$ models. Table~\ref
{table:modFITS} contains estimated parameter values for both training
cases as well as the monitoring observations immediately following the
change-point suggested by the experts. To be precise, the post-change
period is
\begin{enumerate}[(PC)]
\item[(PC)] Observations 18\,501--18\,580.
\end{enumerate}
%
%
\begin{table}[t]
\tablewidth=\textwidth
\tabcolsep=0pt
\caption{Summary of EEG modeling. Standard errors in parenthesis}
\label{table:modFITS}
\begin{tabular*}{\textwidth}{@{\extracolsep{\fill}}lllllld{3.1}@{}}
\hline
Case & $\hat{\phi}_{m,1}$ & $\hat{\phi}_{m,2}$ & $\hat{\phi}_{m,3}$ &
$\hat{\phi}_{m,4}$ & $\hat{\mu}_m$ & \multicolumn{1}{l}{$\hat{\sigma}^2_m$}\\
\hline
(TP1) & 1.66 (0.03)& $-$0.79 (0.06)& $-$0.12 (0.06)& \hphantom{$-$}0.20 (0.03)& $-$207.2 (3.85)&
63.1\\
(TP2) & 1.64 (0.03)& $-$0.74 (0.06)& $-$0.13 (0.06)& \hphantom{$-$}0.18 (0.03)& $-$206.6 (4.90)&
61.9\\
(PC) & 1.46 (0.15)& $-$0.61 (0.27)& \hphantom{$-$}0.20 (0.27)&$-$0.18 (0.16)& $-$194.8 (12.53)&
227.9\\
\hline
\end{tabular*}
\end{table}
All models were fit conditionally on four additional observations in
the respective windows. (E.g., in the case of (TP1), the
$m+p=1004$ observations 15\,997--17\,000 were used for the
estimation.) The tabulated estimates suggest the primary change occurs
in the innovation variance, while the dynamics of the series remains
largely intact.

%
\begin{table}[b]
\tablewidth=\textwidth
\tabcolsep=0pt
\caption{Summary of EEG stopping times and empirical values based on
simulations from the estimated model with 2500 iterations}
\label{table:detections}
\begin{tabular*}{\textwidth}{@{\extracolsep{\fill}}llllllllll@{}}
\hline
& & & & \multicolumn{6}{l}{Simulated empirical values}\\ [-5pt]
& & & & \multicolumn{6}{l}{\hrulefill}\\
& & \multicolumn{2}{l}{Stopping times} & \multicolumn{2}{l}{$95\%$ upper limits} & \multicolumn{2}{l}{Medians} & \multicolumn
{2}{l}{FRR}\\ [-5pt]
& & \multicolumn{8}{l}{\hrulefill}\\
Case & $\gamma$ & Page & \multicolumn{1}{l}{CUSUM}& Page &
\multicolumn{1}{l}{CUSUM}& Page & \multicolumn{1}{l}{CUSUM}& Page &
CUSUM \\
\hline
\multirow{3}{*}{(TP1)} & 0 & 18\,637 & 18\,676&
18\,728& 18\,808 & 18\,623& 18\,643
& 0.0240& 0.0200\\
& 0.25 & 18\,609 & 18\,661 & 18\,718& 18\,798 &18\,614 & 18\,634 & 0.0580& 0.0484\\
& 0.49 & 18\,673 & 18\,691 & 18\,768& 18\,830 & 18\,641& 18\,650 & 0.1344& 0.1296\\
[3pt]
\multirow{3}{*}{(TP2)} & 0 & 18\,580 & 18\,581 &18\,648
& 18\,674 & 18\,576& 18\,588 &
0.0036& 0.0028\\
& 0.25 & 18\,580 & 18\,580 & 18\,626& 18\,657 & 18\,561& 18\,573 &0.0288 & 0.0228\\
& 0.49 & 18\,580 & 18\,580 &18\,633 & 18\,660 &18\,563 & 18\,569 &0.1140 & 0.1104\\
\hline
\end{tabular*}
\end{table}

The proposed testing procedures were applied to the two training sets
at the $\alpha=0.05$ level. Critical values for the CUSUM procedure
were obtained from Horv\'ath \textit{et al.} \cite{hhks04} and critical values for
Page's CUSUM procedure from Fremdt \cite{f12+a}. No
changes were found by the
mean-only procedures $\tau_\mu(m)$ and $\tau_\mu^P(m)$ given in
\eqref{tau_m} and \eqref{tau_m^p} when truncating the tests at
monitoring time point $10m$. The results for the general procedures
$\tau_\bxi(m)$ and $\tau_\bxi^P(m)$ are summarized for three
choices of $\gamma$ in the column labeled ``Stopping Times'' of
Table~\ref{table:detections}. For (TP1), both procedures terminate
within two seconds after the suspected onset of the change, for (TP2)
within one second. Stopping times for (TP1) generally lag behind
stopping times for (TP2). Page's CUSUM detector displays faster
detection for both training periods. 

It should be noted that a sequential procedure does not provide an
estimator for the time of change. In general, it is a difficult problem
to estimate the change-point after a sequential procedure has
terminated because the post-change sample is typically (much) smaller
than the pre-change sample. In the literature, Srivastava and Wu \cite{sw99}
and Wu \cite{w05} have discussed options for this
problem. It would be
worthwhile to follow up on their work elsewhere in the future.


Motivated by the EEG data, several simulations were conducted to
further elaborate on the distribution of the stopping times when a
change occurs only in the innovation variance. The simulations utilized
an $\AR(4)$ model with $\mu= -207$ and $\boldsymbol{\phi}=(1.65,
-0.75, -0.12, 0.18)^\prime$. The innovations were distributed $\Laplace(0, b_0 = 5.6)$ since this closely described the behavior of
the residuals from the EEG training models. Mimicking the two cases
from the EEG application, we used training sizes of $m = 1000$ and
induced changes in the variance by adjusting the scale parameter to
$b_A = 10.7$ at time point 18\,500 (i.e., monitoring time points 500
and 1500 for (TP1) and (TP2), respectively). The choice of scale
parameters imply the difference $\delta^\sigma= 7.21$. Table~\ref
{table:detections} provides simulated empirical confidence limits,
empirical median rejection times and false rejection rates (FRR). The
reported values have been adjusted to fit the time locations observed
in the EEG example. The reported stopping times for the EEG example all
fall within the empirical upper bounds from the simulation study. The
large false rejection rates for $\gamma= 0.49$ display the delay in
convergence to the asymptotic levels suggested by Horv\'{a}th \textit{et al.}
\cite{hhks04} and Fremdt \cite{f12+a} when the sensitivity parameter
is close to
the upper boundary.\looseness=1

\subsection{IBM data}

\begin{figure}[t]

\includegraphics{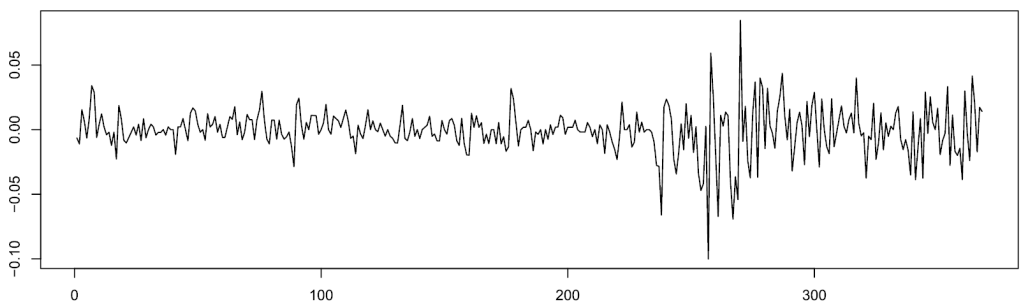}

\caption{Plot of transformed IBM data set.}\vspace*{-6pt}\label{IBMplot}
\end{figure}

The second application is a study of a classic retrospective data set
which has been previously studied for changes in the variance, albeit
in a retrospective setting. The observations are on the IBM common
stock daily closing prices from May 17, 1961 to November 2, 1962. This
is Series B as reported in Box \textit{et al.} \cite{bjr08}. The data set contains 369
observations and has been examined in several retrospective studies
which focused primarily on changes in the variance. Several authors
have detected two change-points. Incl\'{a}n and Tiao \cite{it94} detected
change-points at observations 235 and 279 using their ICSS algorithm,
Baufays and Rasson \cite{br85} proposed 235 and 280, Wichern \textit{et al.} \cite{wmh76}
gave 180 and 235, while Tsay \cite{t88} reported only one
change at
observation 237. As previously suggested in order to stabilize the
variance, the first difference of the log transformed series will be
analyzed. Figure~\ref{IBMplot} displays the corresponding time series
plot. It can be seen that fluctuations appear to be around a constant
level, while amplitudes are larger for roughly the last third of the
observations.
%

To estimate an initial model, the training period is selected to
consist of the first $m=200$ observations. Two competing models were
identified based on AIC and model selection diagnostic plots. The
competing fits are the $\ARMA(2,2)$ and $\AR(4)$ estimated models summarized
in Table~\ref{ibmFITS}, with the AIC value being slightly smaller for
the $\ARMA(2,2)$ model.

\begin{figure}[b]

\includegraphics{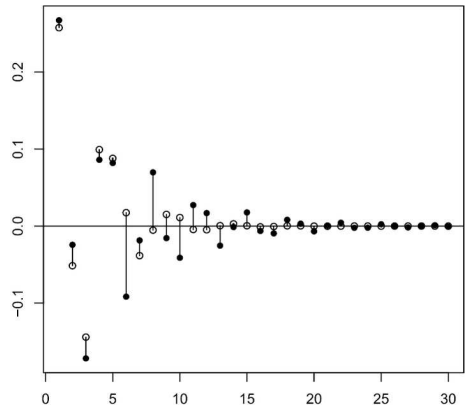}

\caption{Comparing $\MA(\infty)$ coefficients from the
observed $\ARMA(2,2)$ (filled) and $\AR(4)$ (opened)
models.}\label{maINFTY}
\end{figure}

The proposed procedures were applied at the $\alpha=0.05$ level,
utilizing both model fits. Monitoring commences at observation 201. The
mean-only procedures $\tau_\mu(m)$ and $\tau_\mu^P(m)$ do not
detect deviations from a constant level. Table~\ref
{table:IBMdetections} provides the observed values for the general
stopping rules $\tau_\bxi(m)$ and $\tau_\bxi^P(m)$. Depending on
the choice of $\gamma$, both procedures report a change has occurred
at or near observation 238. For comparison purposes, a simulation study
was conducted and is also summarized in Table~\ref
{table:IBMdetections}. The simulations generated training data from the
observed $\ARMA(2,2)$ model. A change was induced at time point
235 to reflect the observed instability in the IBM example. Based on
observations 235--279 (retrospective studies suggest stability over
this period), the best fitting model was a white noise process with
innovation variance given by 0.00135. The empirical measures from the
simulation study are similar when assuming the correct ARMA model
orders, as well as when the $\AR(4)$ is assumed. This highlights
an important feature. Models with nearly identical $\MA(\infty
)$ representations exhibit similar behavior with respect to our
proposed methodology. For our observed $\ARMA(2,2)$ and $\AR(4)$ models, Figure~\ref{maINFTY} displays the differences in the
initial $\MA(\infty)$ coefficients.
%
\begin{table}[b]
\tablewidth=\textwidth
\tabcolsep=0pt
\caption{Summary of IBM modeling. Standard errors in parenthesis}
\label{ibmFITS}
{\fontsize{8}{10}\selectfont{
\begin{tabular*}{\textwidth}{@{\extracolsep{\fill}}lllllllll@{}}
\hline
Model & AIC & $\hat{\phi}_{m,1}$ & $\hat{\phi}_{m,2}$ & $\hat{\phi
}_{m,3}$ & $\hat{\phi}_{m,4}$ & $\hat{\theta}_{m,1}$& $\hat{\theta
}_{m,2}$& $\hat{\sigma}^2_m$\\
\hline
$\ARMA(2,2)$ & $-$1296 & $-$0.40 (0.13)& $-$0.68 (0.11)& \multicolumn{1}{l}{--} & \multicolumn{1}{l}{--} &0.67
(0.12)& 0.76 (0.10)& $8.5\mathrm{e}{-}05$ \\
$\AR(4)$ & $-$1292 & \hphantom{$-$}0.26 (0.07)& $-$0.12 (0.07)& $-$0.10 (0.07)& 0.16 (0.07)&
\multicolumn{1}{l}{--} & \multicolumn{1}{l}{--} & $8.7\mathrm{e}{-}05$ \\
\hline
\end{tabular*}}}%
\end{table}
%
\begin{table}[t]
\tablewidth=\textwidth
\tabcolsep=0pt
\caption{Summary of IBM stopping times and empirical values based on
simulations from the estimated model with 2500 iterations}
\label{table:IBMdetections}
\begin{tabular*}{\textwidth}{@{\extracolsep{\fill}}llllllllll@{}}
\hline
& & & & \multicolumn{6}{l}{Simulated empirical values}\\ [-5pt]
& & & & \multicolumn{6}{l}{\hrulefill}\\
& & \multicolumn{2}{l}{Stopping times} & \multicolumn{2}{l}{$95\%$ upper limits} & \multicolumn{2}{l}{Medians} & \multicolumn
{2}{l}{FRR}\\ [-5pt]
& & \multicolumn{8}{l}{\hrulefill}\\
Case & $\gamma$ & Page & \multicolumn{1}{l}{CUSUM}& Page &
\multicolumn{1}{l}{CUSUM}& Page & \multicolumn{1}{l}{CUSUM}& Page &
CUSUM \\
\hline
\multirow{3}{*}{$\ARMA(2,2)$} & 0 & 238 & 239& 244 & 244 &238 & 238 &0.0024 &0.0024 \\
& 0.25 & 238 & 238 & 242 & 242 &237 & 237 & 0.0228& 0.0216 \\
& 0.49 & 238 & 238 & 241 & 242 & 236& 237 & 0.1008& 0.1072\\
[3pt]
\multirow{3}{*}{$\AR(4)$} & 0 & 239 &242
& 244 & 244 &238 & 239
&0.0004&0.0004 \\
& 0.25 & 238 & 238 &242 & 243 &237 & 237 &0.0120& 0.0100\\
& 0.49 & 238 & 238 & 241 & 242 &237 & 237 & 0.0848& 0.0904\\ \hline
\end{tabular*}
\end{table}

%

\section{Proofs}\label{sec5}

\subsection{Preliminaries}

The following auxiliary result will be used frequently. It establishes
the behavior of the coefficients $\pi_j(\bv)$ and $\psi_j(\bu)$ in
\eqref{recip} if instead of the true parameter vectors $\bphi$ and
$\btheta$, generic elements $\bv\in\R^p$ and $\bu\in\R^q$ in
their vicinity are used for the respective power series expansions. Let
$|\cdot|$ denote the maximum norm of vectors.

\begin{prop}\label{prop1}
Let $(Y_t\dvtx  t\in\mathbb{Z})$ follow the ARMA equations \eqref
{ARMA} so that \eqref{inv-cau} holds. Let $\bv\in\R^p$ and $\bu
,\bu_1,\bu_2\in\R^q$. Then there are $\eps> 0, c\in(0,1)$ and
$K>0$ such that, for all $j\geq0$,
\begin{enumerate}[(b)]
\item[(a)] $|\pi_j(\bv)|\leq K c^j, \mbox{ if }|\bv-
\bphi|\leq \eps$;

\item[(b)] $|\psi_j(\bu)|\leq K c^j, \mbox{ if }|\bu-
\btheta |\leq\eps$;

\item[(c)] $|\psi_j(\bu_1) - \psi_j(
\bu_2)|\leq K |\bu_1-\bu _2|jc^{j-1},
\mbox{ if }|\bu_1 - \btheta|\leq\eps\mbox{ and }|\bu _2
- \btheta|\leq\eps$.
\end{enumerate}
\end{prop}

\begin{pf}
The proof of these statements can be found in Bai \cite{b93}.
\end{pf}

Throughout the proofs, let $\lambda_t$ denote the difference between
the residuals $\epsh_t$ and the innovations $\eps_t$ if the null
hypothesis $H_0$ is valid. Since none of the parameters is subject to
change, it holds then that
%
\begin{equation}
\label{lambda} \lambda_t = \lambda_t \bigl(\sqrt{m}[
\hat{\boldsymbol{\theta }}_{m}- \btheta_0],\sqrt {m}[\hat{
\boldsymbol{\phi}}_{m}- \bphi_0],\sqrt{m}[\hat{\mu
}_{m}- \mu_0] \bigr),
\end{equation}
where
\begin{eqnarray*}
\lambda_t(\bu,\bv,w) = \zeta_t(\bu)+\frac{\beta_t(\bu,\bv
)}{\sqrt{m}}
+ \frac{\rho_t(\bu,\bv,w)}{\sqrt{m}}.
\end{eqnarray*}
To define the quantities on the right-hand side of the latter equality,
let first $\bu^* = \btheta+ \bu/\sqrt{m}$ and $\bv^* = \bphi+ \bv
/\sqrt{m}$, and set $u_0^*=0$. Then
\begin{eqnarray*}
\zeta_t({\bu}) &=&-\sum_{j=1}^q
\Biggl(\sum_{\ell=1}^j\psi_{t-1+\ell}
\bigl({\bu }^*\bigr)u_{j-\ell}^* \Biggr)\eps_{1-j},
\\
 \beta_t({\bu},{\bv}) &=&-\sum_{j=1}^pv_j
\sum_{\ell=0}^{t-1}\psi_\ell\bigl({\bu
}^*\bigr)X_{t-j-\ell}-\sum_{j=1}^qu_j
\sum_{\ell=0}^{t-1}\psi_\ell\bigl({\bu
}^*\bigr)\eps_{t-j-\ell},
\\
\rho_t({\bu},{\bv},w) &=&- \Biggl(1-\sum
_{j=1}^pv_j^* \Biggr)w \sum
_{\ell=0}^{t-1}\psi_\ell \bigl({\bu}^*\bigr).
\end{eqnarray*}
This is the decomposition given in Yu \cite{y07}, which
is useful to derive
the limit distributions of the various test procedures under the null
hypotheses as given in Theorems \ref{Th1} and \ref{Th2}. This was
done for the CUSUM-type procedure in Dienes and Aue \cite{da12+}, but the
same approach works also for the procedure based on Page's CUSUM using
the work of Fremdt \cite{f12+a,f12+b}. 

To prove the new results on the asymptotic delay time distribution of
the stopping times one may modify methodology developed in Aue and Horv{\'a}th \cite{ah04}:
It is subsequently shown that sequences $N=N(m,x)$
can be found such that, for the stopping time $\tau$ with
corresponding detector $D(m,k)$, it holds that
\begin{eqnarray*}
P(\tau>N) = P \biggl(\max_{1\leq k\leq N}\frac{\hat D(m,k)}{g_\gamma(m,k)} \leq c
\biggr)
\end{eqnarray*}
converges to the appropriate limit distribution. The standardizations
for $\tau$ in the various theorems are then implied by the definition
of $N$. The next section contains the verification for the mean break case.

\subsection{Proofs of the results in Section~\texorpdfstring{\protect\ref{ssec1}}{3.2}}\label{Pssec1}

For the mean break case, changes in the second order parameters $\bphi
$, $\btheta$ and $\sigma^2$ are precluded. To determine the effect of
the mean break on the differences $\epsh_t-\eps_t$, one consequently
needs to check only the terms including the $\mu_{t-j}$. It can be
seen from \eqref{lambda} that these terms only enter through $\rho
_t$. To determine the drift induced by the change in mean under
$H_A^\mu$, a similar decomposition to \eqref{lambda} is needed.
Following equation (14) in Yu \cite{y07}, it follows that
%
\begin{eqnarray}\label{Yu14}
\epsh_t - \eps_t &=& \tilde\zeta_t(\hat{\boldsymbol{\theta}}_{m})
+ \frac{\tilde
\beta_t(\hat{\boldsymbol{\theta}}_{m}
,\hat{\boldsymbol{\phi}}_{m})}{\sqrt{m}}\nonumber \\[-8pt]\\[-8pt]
&&{}- \sum_{\ell= 0}^{t-1}
\psi_\ell(\hat{\boldsymbol{\theta }}_{m}) \Biggl[(\hat{
\mu}_{m}- \mu _{t-k}) - \sum_{j = 1}^p
\hat{\phi}_{m,j}(\hat{\mu}_{m}-\mu _{t-j-\ell}) \Biggr],\nonumber
\end{eqnarray}
where $\tilde\zeta_t(\hat{\boldsymbol{\theta}}_{m})=\zeta_t(\sqrt {m}[\hat{\boldsymbol{\theta}}_{m}
-\btheta_0])$ and $\tilde\beta_t(\hat{\boldsymbol{\theta
}}_{m},\hat{\boldsymbol{\phi}}_{m})=\beta_t(\sqrt {m}[\hat{\boldsymbol{\theta}}_{m}-\btheta_0],\sqrt{m}[\hat
{\boldsymbol{\phi}}_{m}-\bphi_0])$ are
respectively the terms of initialization effects and the partial sums
of centered observations and innovations. To derive \eqref{Yu14}, one
uses the recursiveness of the difference $\epsh_t-\eps_t$ and the
invertibility of the underlying ARMA process.
Now, as under $H_A^\mu$ a change occurs only in $\mu_t$ for $t\geq
m+{k^{*}}$, it suffices to investigate the term
\begin{eqnarray*}
&&- \sum_{\ell= 0}^{t-1} \psi_\ell(
\hat{\boldsymbol{\theta}}_{m}) \Biggl[ (\hat{\mu }_{m}-
\mu_{t-k}) - \sum_{j = 1}^p \hat{
\phi}_{m,j}(\hat{\mu}_{m}- \mu_{t-j-\ell}) \Biggr]
\\
&&\quad =- \Biggl(1 -\sum_{j = 1}^p \hat{
\phi}_{m,j} \Biggr) (\hat{\mu }_{m}-\mu_0)\sum
_{\ell= 0}^{t-1}\psi_\ell(\hat{
\boldsymbol{\theta}}_{m}) +\delta_{m}^\mu\sum
_{\ell= 0}^{t-1}\psi_\ell(\hat{
\boldsymbol {\theta}}_{m}) \Biggl[I_{t,0,\ell} - \sum
_{j = 1}^p \hat{\phi}_{m,j}I_{t,j,\ell}
\Biggr]
\\
&&\quad =\frac{\tilde\rho_t(\hat{\boldsymbol{\theta}}_{m},\hat
{\boldsymbol{\phi}}_{m},\hat{\mu}_{m})}{\sqrt{m}}+\Lambda _{t-m-{k^{*}}}^\mu,
\end{eqnarray*}
where $\tilde\rho_t(\hat{\boldsymbol{\theta}}_{m},\hat
{\boldsymbol{\phi}}_{m},\hat{\mu}_{m})=\rho_t(\sqrt {m}[\hat{\boldsymbol{\theta}}_{m}-\btheta_0],\sqrt{m}[\hat
{\boldsymbol{\phi}}_{m}-\bphi_0],\sqrt{m}[\hat{\mu}_{m}
-\mu_0])$ and $I_{t,j,\ell}$ is short for $I_{\{t-j-\ell\geq{k^{*}}+
m\}}(t,j,\ell)$. Here, $I_A$ denotes the indicator function of a set
$A$. Letting $t\geq m+{k^{*}}$ and $s = t-m-{k^{*}}$, the drift term
can be
written as\vspace*{-1pt}
%
\begin{eqnarray}\label
{Lambda^mu}
\Lambda^\mu_s &=& \delta_m^\mu\sum
_{\ell= 0}^{t-1}\psi_\ell(\hat{
\boldsymbol {\theta}}_{m}) \Biggl[I_{t,0,\ell} - \sum
_{j = 1}^p \hat{\phi}_{m,j}I_{t,j,\ell}
\Biggr] \nonumber
\\[-8pt]\\[-8pt]
&=&\cases{ 0,&\quad  $s<0$,
\cr
\delta_m^\mu\displaystyle \sum
_{\ell= 0}^{s}\psi_{s-\ell}(\hat{\boldsymbol {
\theta}}_{m}) \Biggl(1 - \displaystyle \sum_{j = 1}^\ell
\hat{\phi}_{m,j} \Biggr),& \quad $0\leq s<p$,
\cr
\delta_m^\mu
\Biggl[ \Biggl(1 - \displaystyle \sum_{j = 1}^p \hat{
\phi}_{m,j} \Biggr)\sum_{\ell= 0}^{s-p}
\psi_\ell(\hat{\boldsymbol{\theta}}_{m}) + \displaystyle \sum
_{\ell= 0}^{p-1}\psi_{s-\ell}(\hat{\boldsymbol{
\theta }}_{m}) \Biggl(1 - \displaystyle \sum_{j =
1}^\ell
\hat{\phi}_{m,j} \Biggr) \Biggr], &\quad  $s\geq p$. }
\nonumber
\end{eqnarray}
Note that the drift has been rescaled, so that $s<0$ indicates that the
change has not yet occurred. The further distinction into the cases
$0\leq s<p$ and $s\geq p$ takes into account the autoregressive order.
It follows that\vspace*{-1pt}
%
\begin{equation}\label{epsdiff}
\epsh_t - \eps_t = \lambda_t +
\Lambda^\mu_{t-m-{k^{*}}},
\end{equation}
with $\lambda_t$ from \eqref{lambda}. To prove the theorems of
Section~\ref{ssec1}, it remains to analyze partial sums of the $ \epsh
_t - \eps_t$ and compare them to the growth of the threshold $g_\gamma(m,k)$.

\begin{pf*}{Proof of Theorem~\ref{AD}}
Let $k \geq{k^{*}}$ and $M = k-{k^{*}}$. Utilizing $\lambda_t$ from
\eqref
{lambda} and display \eqref{epsdiff}, it follows that\vspace*{-1pt}
\begin{eqnarray*}
\sum_{t= m+1}^{m+k} (\epsh_t-
\eps_t) = \sum_{t=m+1}^{m+k}
\bigl(\lambda_t + \Lambda^\mu_{t-m-{k^{*}}}\bigr) = \sum
_{t = m+1}^{m+k} \lambda_t + \sum
_{s = 0}^M \Lambda^\mu_s.
\end{eqnarray*}
The first term on the right-hand side can be treated as under the null
hypothesis, see Dienes and Aue \cite{da12+}. The drift of the cumulative sum
procedure can be determined as follows. First, for $M<p$, \eqref
{Lambda^mu} implies directly that\vspace*{-1.5pt}
\[
\sum_{s = 0}^M\Lambda^\mu_s
= \delta_m^\mu\sum_{s = 0}^M
\sum_{\ell= 0}^{s}\psi_{s-\ell
}(\hat{
\boldsymbol{\theta}}_{m}) \Biggl(1 - \sum_{j = 1}^\ell
\hat {\phi}_{m,j} \Biggr).
\]
Second, for $M\geq p$, another application of \eqref{Lambda^mu} using
the cases for $0\leq s<p$ and $p\leq s\leq M$ to split up the sum and
subsequently combining the terms involving the incomplete sums
$1-\sum_{j=1}^\ell\hat\phi_{m,j}$ of estimated autoregressive coefficients, yields\vspace*{-1.5pt}
\begin{eqnarray*}
\sum_{s = 0}^M\Lambda^\mu_s
&=& \delta_m^\mu
\Biggl[ \Biggl(1 - \sum_{j = 1}^p \hat{\phi
}_{m,j} \Biggr)\sum_{\ell= 0}^{M-p}
\psi_\ell(\hat{\boldsymbol{\theta }}_{m}) \bigl[(M-p+1) -
\ell \bigr]
\\
&&\hphantom{\delta_m^\mu
\Biggl[}{} + \sum_{s = 0}^{p-1} \Biggl(1 - \sum
_{j = 1}^s \hat{\phi }_{m,j} \Biggr)
\sum_{\ell= 0}^{M-s}\psi_{\ell}(\hat{
\boldsymbol{\theta }}_{m}) \Biggr]
\\
 &=& \hat\Delta_m^\mu(M-p+1) + \delta_m^\mu
\bigl[A_1(M) - A_2(M) + A_3(M) \bigr],
\end{eqnarray*}
where $\hat\Delta_m^\mu=\delta_m^\mu\hat\phi_{m}(1)/\hat\theta
_m(1)$ with $\hat\phi_m(1)=1-\hat{\phi}_{m,1}z-\cdots-\hat{\phi
}_{m,p}z^p$ and
$\hat{\theta}_{m}(1)=1+\hat{\theta}_{m,1}z+\cdots+\hat{\theta
}_{m,q}z^q$, and
\begin{eqnarray*}
A_1(M) &=& (M-p+1) \Biggl(1 - \sum_{j = 1}^p
\hat{\phi}_{m,j} \Biggr) \sum_{\ell= M-p+1}^{\infty}
\psi_\ell(\hat{\boldsymbol{\theta }}_{m}),
\\
A_2(M) &=& \Biggl(1 - \sum_{j = 1}^p
\hat{\phi}_{m,j} \Biggr)\sum_{\ell=
0}^{M-p}
\ell\psi_\ell(\hat{\boldsymbol{\theta}}_{m}),
\\
A_3(M) &=& \sum_{s = 0}^{p-1}
\Biggl(1 - \sum_{j = 1}^s \hat{\phi
}_{m,j} \Biggr)\sum_{\ell= 0}^{M-s}
\psi_{\ell}(\hat{\boldsymbol {\theta}}_{m}).
\end{eqnarray*}
It is clear that $\hat\Delta_m^\mu$ will be close to its
deterministic equivalent $\Delta_m^\mu$ if $m$ is large. The terms
$A_1(M)$, $A_2(M)$ and $A_3(M)$ are stochastically bounded, so that
Proposition~\ref{prop1} implies that, as $m\to\infty$,
\begin{eqnarray*}
\biggl(\frac{N}{m} \biggr)^{\g-\oh}\max_{{k^{*}}\leq k\leq N}
\frac
{\delta_m^\mu A_i(k-{k^{*}}
)}{g_\gamma(m,k)} = \mathrm{o}_P(1),\qquad  i = 1,2,3.
\end{eqnarray*}
If the sequence $N = N(m,x)$ given in Fremdt \cite{f12+a}
is used as the
upper bound for the maximum. The rest of the proof of part (a) of the
theorem follows now analogously to the proof of Theorem~2.2 in
Fremdt \cite{f12+a}.

Part (b) can be verified by an extension of the proof in Aue and Horv\'
ath \cite{ah04}, relaxing their assumption on the order of the change-point
to the requirement of part (a) in Assumption~\ref{As1}. This can be
done with the sharper estimates developed in Fremdt \cite
{f12+a}. Further
details are omitted to conserve space.
\end{pf*}

\begin{pf*}{Proof of Theorem~\ref{ADmu^2}}
To investigate the behavior of the general detectors under $H_A^\mu$,
the previous proof needs to be adjusted for the squared residuals. From
\eqref{epsdiff} it follows that, for $t \geq m+{k^{*}}$,
%
\begin{eqnarray}
\label{ressq} \epsh_t^2 - \eps_{t}^2
= \lambda_t^2 + 2\lambda_t
\eps_{t} + \bigl(\Lambda^\mu_{t-m-{k^{*}}}
\bigr)^2 + 2\Lambda^\mu_{t-m-{k^{*}}}(\eps
_{t} + \lambda_t ).
\end{eqnarray}
The first two terms on the right-hand side can again be treated as
under the null hypothesis. The relevant drift term for the sequential
procedures consists then of the partial sums of $(\Lambda^\mu
_{t-m-{k^{*}}})^2$ and $2\Lambda^\mu_{t-m-{k^{*}}}(\eps_{t} +
\lambda_t
)$, of which the latter will be negligible. To verify this claim,
observe first that
\begin{eqnarray*}
\max_{{k^{*}}\leq k<\infty} \frac{1}k \sum
_{t = 0}^k |\eps_t| =
\mathrm{O}_P(1)
\end{eqnarray*}
since, on account of the strong law of large numbers, $\frac{1}k \sum_{t = 0}^k |\eps_t|$ converges almost surely as $k\to\infty$. Because
\begin{eqnarray*}
\max_{{k^{*}}\leq k\leq N} \biggl(\frac{N}{m} \biggr)^{\g-\oh}
\frac
{k}{g_\gamma(m,k)} = \mathrm{o} ({1} )\qquad  (m\to \infty),
\end{eqnarray*}
it follows from Proposition~\ref{prop1} that
\begin{eqnarray*}
\biggl(\frac{N}{m} \biggr)^{\g-\oh}\max_{{k^{*}}\leq k\leq N}
\sum_{t = m+{k^{*}}}^{m+k}\frac{\Lambda
_{t-m-{k^{*}}}^\mu\eps_{t}}{g_\gamma(m,k)} =
\mathrm{o}_P(1) \qquad (m\to\infty).
\end{eqnarray*}
Utilizing the definition of $\Lambda^\mu_s$ in \eqref{Lambda^mu} and
another application of Proposition~\ref{prop1} in combination with
Lemmas 6.1--6.3 of Dienes and Aue \cite{da12+} yield also that
\begin{eqnarray*}
\biggl(\frac{N}{m} \biggr)^{\g-\oh}\max_{{k^{*}}\leq k\leq N}
\sum_{t =m+ {k^{*}}}^{m+k}\frac{\Lambda
_{t-m-{k^{*}}}^\mu\lambda_t}{g_\gamma(m,k)} =
\mathrm{o}_P(1) \qquad (m\to\infty).
\end{eqnarray*}
It therefore remains to extract the dominating term from the partial
sums of $(\Lambda_{t-m-{k^{*}}}^\mu)^2$. To facilitate notation, the
abbreviations $\psi_\ell= \psi_\ell(\hat{\theta}_m)$, $\hat\phi
_m^{(\ell)}(z)=1-\hat{\phi}_{m,1}z-\cdots-\hat{\phi}_{m,\ell
}z^\ell$, $\ell
=1,\ldots,p-1$, and $k' = {k^{*}}+m+p$ are used. Then,
\begin{eqnarray*}
\sum_{t = k'}^{m+k} \bigl(
\Lambda_{t-m-{k^{*}}}^\mu \bigr)^2 &=& \bigl(\delta
_m^\mu\bigr)^2 \Biggl[  \hat
\phi_m^2(1)
\sum_{t = k'}^{m+k}
\Biggl(\sum_{\ell= 0}^{t-k'}\psi_\ell
\Biggr)^2 +\sum_{t = k'}^{m+k}
\Biggl(\sum_{\ell= 0}^{p-1}\hat
\phi_m^{(\ell
)}(1)\psi_{t- m -{k^{*}}-\ell} \Biggr)^2
\\
&&\hphantom{\bigl(\delta
_m^\mu\bigr)^2 \Biggl[}{}-2\hat\phi_m^2(1)
\sum
_{t = k'}^{m+k} \Biggl(\sum_{\ell= 0}^{t-k'}
\psi_\ell \Biggr) \Biggl(\sum_{\ell' = 0}^{p-1}
\hat\phi_m^{(\ell')}(1)\psi _{t-m-{k^{*}}-\ell'} 
\Biggr)
\Biggr].
\end{eqnarray*}
Similar arguments to those used in the proof of Theorem~\ref{AD} yield
that only the first term needs to be investigated. Since
\begin{eqnarray*}
\sum_{s = 0}^{t} \Biggl(\sum
_{\ell= 0}^{s}\psi_\ell \Biggr)^2
&=& (t+1) \Biggl(\sum_{\ell= 0}^{\infty}
\psi_\ell \Biggr)^2 - 2 \Biggl(\sum
_{\ell= 0}^{\infty}\psi_\ell \Biggr)\sum
_{s = 0}^{t} \Biggl(\sum_{\ell= s+1}^{\infty}
\psi_\ell \Biggr) + \sum_{s = 0}^{t}
\Biggl(\sum_{\ell= s+1}^{\infty}\psi_\ell
\Biggr)^2
\\
&=& (t+1) \Biggl(\sum_{\ell= 0}^{\infty}
\psi_\ell \Biggr)^2 - 2 \Biggl(\sum
_{\ell= 0}^{\infty}\psi_\ell \Biggr) \Biggl[\sum
_{s =
1}^{t}s\psi_s + (t+1)\sum
_{s = t+1}^{\infty}\psi_s \Biggr]
\\
&&{} + \sum_{s = 0}^{t} \Biggl(\sum
_{\ell= s+1}^{\infty}\psi_\ell \Biggr)^2,
\end{eqnarray*}
following the arguments of the proof of Theorem~\ref{AD} implies that
$(t+1)\hat{\theta}_{m}^{-2}(1)$ is the dominating term in this
expression. The
rest follows analogously to the proof of Theorem~2.2 in Fremdt \cite{f12+a}.
\end{pf*}

\subsection{Proofs of the results in Section~\texorpdfstring{\protect\ref{ssec2}}{3.3}}\label{Pssec2}

Denote by $(z_t\dvtx  t\in\mathbb{Z})$ the sequence of independent,
identically distributed and standardized random variables given by the
requirement $\eps_t=\sigma_tz_t$ for all $t\in\mathbb{Z}$.
Therefore changes in scale do not affect the $z_t$'s. In the following,
the subscript 0 in the quantities $y_{0,t}$, $x_{0,t}$ and $\eps
_{0,t}$ will indicate that the corresponding random variables are
generated according to the null parameter vector $\bxi_0$.

Precluding a break in the mean, autoregressive and moving average
parameters and only allowing breaks in the scale parameter, leads to
the decomposition
%
\begin{eqnarray}
\label{ressq-var} \epsh_t^2 - \eps_{0,t}^2
= \lambda_t^2 + 2\lambda_t
\eps_{0,t} + \bigl(\Lambda^{\sigma}_{t-m-k^*}
\bigr)^2 + 2\Lambda^{\sigma}_{t-m-k^*}(\eps
_{0,t}+\lambda_t),
\end{eqnarray}
for $t\geq m+k^*$, which is analogous to \eqref{ressq}. Now $\Lambda
_{t-m-k^*}^\sigma$ can be further decomposed into 
%
\begin{eqnarray*}
\Lambda_{t-m-k^*}^\sigma= \delta_m^\sigma(z_{t}
+ B_t).
\end{eqnarray*}
%
Using that, for $t\geq m+k^*$,
\begin{eqnarray*}
x_t - x_{0,t} = \delta_m^\sigma
\Biggl(\sum_{k = 0}^{t-m-{k^{*}}} \pi_k(
\bphi _A) \Biggl[z_{t-k} + \sum
_{j = 1}^{\min(q,t-m-{k^{*}}-k)} \theta _{0,j}z_{t-j-k}
\Biggr] \Biggr),
\end{eqnarray*}
and setting again $s=t-m-k^*$ gives
\begin{eqnarray*}
B_t &= &\sum_{j=1}^q(
\theta_{0,j}-\hat{\theta}_{m,j})\sum
_{k=0}^{s-q}\psi _kz_{t-j-k} +
\sum_{k=1}^{q-1}\sum
_{j=1}^{k}\psi_{s-k}(\theta
_{0,j}-\hat{\theta}_{m,j})z_{m+{k^{*}}-j+k}
\\
 &&{}+ \sum_{j=1}^p(
\phi_{0,j}-\hat{\phi}_{m,j})\sum_{k=0}^{s-j}
\psi _k\sum_{n=0}^{s-j-k}
\pi_n(\bphi_0)z_{t-j-k-n}
\\
 &&{} +\sum_{j=1}^p(
\phi_{0,j}-\hat{\phi}_{m,j})\sum_{\ell
=1}^{q}
\theta _{0,j}\sum_{k=0}^{s-j}\sum
_{n=0}^{s-j-k-q}\psi_k
\pi_n(\bphi _0)z_{t-j-k-\ell-n}
\\
 &&{} +\sum_{j=1}^p(
\phi_{0,j}-\hat{\phi}_{m,j})\sum_{k=0}^{s-j}
\psi _k\sum_{n=1}^{q-1}
\pi_{s-j-k-n}(\bphi_0)\sum_{\ell=1}^{n}
\theta_{0,\ell
}z_{m+{k^{*}}+n-\ell}
\\
 &=&B_{1,t}+\cdots+B_{5,t},
\end{eqnarray*}
where $\psi_k=\pi_k(\bphi_0)=0$ for $k<0$. The following lemma
identifies the dominating term in the partial sums of $\hat\eps
_t^2-\eps_{0,t}^2$.
%

\begin{Le}\label{LE}
Under the assumptions of Theorem~\ref{AD-sigma},
\[
\biggl(\frac{N}{m} \biggr)^{\g-\oh} \max_{{k^{*}}\leq k\leq N}
\frac{1}{g_\gamma(m,k)} \Biggl|\sum_{t=m+k^*}^{m+k}\bigl(
\hat \eps_t^2-\eps_{0,t}^2\bigr)-
\Delta_m^\sigma\sum_{t=m+k^*}^{m+k}z_t^2
\Biggr| =\mathrm{o}_P(1)
\]
as $m\to\infty$.
\end{Le}

\begin{pf}
It suffices to examine the quantities on the right-hand side of \eqref
{ressq-var}. Notice first that $\lambda_t^2+2\Delta_t\eps_{0,t}$
contains only terms related to the behavior under the null hypothesis.
For the next two terms on the right-hand side of \eqref{ressq-var}, write
%
\begin{eqnarray}
\label{decomp} &&\bigl(\Lambda_s^\sigma\bigr)^2+2
\Lambda_s^\sigma(\eps_{0,t}+\lambda_t)\nonumber\\[-8pt]\\[-8pt]
&&\quad =\bigl[\bigl(\delta_m^\sigma\bigr)^2+2
\sigma_0\delta_m^\sigma\bigr]z_t^2
+2\delta_m^\sigma\bigl(\delta_m^\sigma+
\sigma_0\bigr)z_tB_t +\bigl(
\delta_m^\sigma\bigr)^2B_t^2
+2\lambda_t\Lambda_s^\sigma.\nonumber
\end{eqnarray}
The first term is the dominating term. Since $\Delta_m^\sigma=(\delta
_m^\sigma)^2+2\sigma_0\delta_m^2$, the assertion of the lemma will
follow if the remaining terms can be shown to be negligible. For the
second term notice that
\[
\biggl(\frac{N}{m} \biggr)^{\g-\oh} \max_{{k^{*}}\leq k\leq N}
\sum_{t=m+k^*}^{m+k}\frac{z_tB_t}{g_\gamma(m,k)} \leq
\biggl(\frac{N}{m} \biggr)^{\g-\oh} \sum
_{t=m+k^*}^{m+N}\frac{|z_t||B_t|}{g_\gamma(m,k^*)} =\mathrm{o}_P(1),
\]
since $z_t$ and $B_t$ are independent and $\sqrt{m}E[B_t]<\infty$,
following the arguments used in Dienes and Aue \cite{da12+}. For the third
term in \eqref{ressq-var}, observe that
\[
\biggl(\frac{N}{m} \biggr)^{\g-\oh} \max_{{k^{*}}\leq k\leq N}
\sum_{t=m+k^*}^{m+k}\frac{B_t^2}{g_\gamma(m,k)} \leq
\biggl(\frac{N}{m} \biggr)^{\g-\oh}\sum_{t=m+k^*}^{m+k}
\sum_{\ell=1}^5 \frac{B_{\ell,t}^2}{g_\gamma(m,k^*)}.
\]
The proof is only detailed for $\ell=1$, since all other terms can be
handled in a similar fashion. For this case,
\begin{eqnarray*}
&&\biggl(\frac{N}{m} \biggr)^{\g-\oh} \frac{1}{g_\gamma(m,k^*)}\sum
_{t=m+k^*}^{m+N} \Biggl(\sum_{j=1}^q(
\theta_{0,j}-\hat{\theta}_{m,j})\sum
_{k=0}^{s-q}\psi _kz_{t-j-k}
\Biggr)^2
\\
&&\quad  \leq \biggl(\frac{N}{m} \biggr)^{\g-\oh} \frac{q}{g_\gamma
(m,k^*)}
\sum_{t=m+k^*}^{m+N}\sum
_{j=1}^q(\theta_{0,j}-\hat{
\theta}_{m,j})^2 \Biggl(\sum_{k=0}^{s-q}
\psi_kz_{t-j-k} \Biggr)^2.
\end{eqnarray*}
Now the arguments of Lemma~5.2 in Dienes and Aue \cite{da12+} apply and
yield the $\mathrm{o}_P(1)$ rate. For the last term in \eqref{ressq-var} there
is nothing to show, since
\[
\Lambda_s^\sigma\lambda_t =
\delta_m^\sigma(\lambda_tz_t+
\lambda_tB_t) \leq\delta_m^\sigma
\bigl(\lambda_tz_t+\lambda_t^2+B_t^2
\bigr)
\]
and all these terms have already been shown to be negligible. The proof
is complete.
\end{pf}

\begin{pf*}{Proof of Theorem~\ref{AD-sigma}}
The relevant drift term has been identified in Lemma~\ref{LE}.
Noticing that the law of the iterated logarithm implies that, for all
$\delta\in(0,1/2)$ and as $m\to\infty$,
\begin{eqnarray*}
&&\biggl(\frac{N}{m} \biggr)^{\g-\oh} \max_{k^*\leq k\leq N}
\frac
{1}{g_\gamma
(m,k)}\sum_{t=m+k^*}^{m+k}
\bigl(z_t^2-1\bigr)
\\
&&\quad =\mathrm{O}_P(1)\frac{m^{1-\gamma}}{N^{1/2-\gamma}}\max_{1\leq
k\leq N-k^*}
\frac{k^{1/2-\gamma+\delta}}{(m+k^*+m)^{1-\gamma}} =\mathrm{o}_P(1),
\end{eqnarray*}
the assertion of the theorem follows.
\end{pf*}

\section*{Acknowledgements}
This research was partially supported by Collaborative Research Center ``Statistical modeling of nonlinear dynamic processes''
(SFB 823) of the German Research Foundation (DFG), NSF grants DMS 0905400, DMS
1209226 and DMS 1305858, and DFG grant STE 306/22-1.


%

\printhistory

\end{document}